\documentclass{article}

\usepackage{arxiv}

\usepackage[utf8]{inputenc} 
\usepackage[T1]{fontenc}    
\usepackage{hyperref}       
\usepackage{url}            
\usepackage{booktabs}       
\usepackage{amsfonts}       
\usepackage{nicefrac}       
\usepackage{microtype}      
\usepackage{lipsum}

\usepackage{subcaption}
\usepackage{amsmath,amssymb,color,amsthm}
\usepackage{epsfig}
\usepackage{algorithm2e}
\usepackage{booktabs}
\usepackage{units}

\newtheorem{definition}{Definition}

\def\R{\hbox{$\mathbb R$}}

\def\N{\hbox{$\mathbb N$}}
\def\X{\hbox{$\mathbb X$}}

\title{A Ray Tracing Technique for the Navigation on a Non-convex Pareto Front}

\author{
  Dimitri Nowak \\
  Optimization Department\\
  Fraunhofer Institute for Industrial Mathematics ITWM\\
  Fraunhofer-Platz 1, 67663 Kaiserslautern, Germany \\
  \texttt{dimitri.nowak@itwm.fraunhofer.de} \\
   \And
 Karl-Heinz K{\"u}fer \\
  Optimization Department\\
  Fraunhofer Institute for Industrial Mathematics ITWM\\
  Fraunhofer-Platz 1, 67663 Kaiserslautern, Germany \\
}

\begin{document}
\maketitle

\begin{abstract}
A new interactive approach to navigate on approximations of in general non-convex but connected Pareto fronts is introduced. Given a finite number of precalculated representative Pareto-efficient solutions, an adapted Delaunay triangulation is generated. Based on interpolation and ray tracing techniques, real time navigation in the vicinity of Pareto-optimal solutions is made possible.
\end{abstract}



\section{Introduction}
The industry offers a wide range of optimization problems with conflicting objectives. More often than not these problems do not have one universal solution. Solutions that perform best in one objective can be improved in another objective by losing some performance in the one objective. A trade-off has to be made dependent on the decision maker. The main issue is, though, that the decision maker often does not know what {trade-off} may be best for him. At the Fraunhofer Institute ITWM, it is standard practice to provide decision makers with alternative solutions and to support with navigation tools that allow an exploration of the entire solution space \cite{Monz2008}.

In the literature \cite{Hwang1979}, different decision support methods are being suggested. In general, they are subdivided into four categories. {A priori methods are based on the decision makers knowledge \cite{Harrington1965}, e.g., how to weight his objectives in order to incorporate them into an appropriate optimization problem.} The problem is solved and the decision maker is presented with the best outcome. A posteriori methods compute some {Pareto optimal} solutions first. Then, all are presented to the decision maker {who chooses} {the most preferred}. No-preference-methods do not include any input from the decision maker. And, finally, interactive methods allow the decision maker to change his preference during an interactive decision process by exploring the entire space of Pareto optimal solutions. Some methods are {designed} to explore a convex Pareto front. They usually perform faster than those applicable to nonconvex problems. However, nonconvex methods can also be applied to convex problems, whereas, convex methods can not process nonconvex information properly. {Most recent review on different types of navigation strategies and some application areas can be found in \cite{Allmendinger2017}}.

The method presented in this paper is an interactive method that is specifically applied to nonconvex multicriteria optimization problems. It can be considered as a three stage process. In the first stage, Pareto optimal solutions have to be computed by algorithms that yield a good representation of the original Pareto front. In the second stage, an approximation of the Pareto front is constructed and serves as a surrogate model of a presumably computationally expensive original problem. Finally, in the third stage, the surrogate model is used to explore an approximation of the actual Pareto front in terms of an interactive behavior of the decision maker such as selection and restriction of both objective and parameter values.

The first stage of the process involves a combination of Sandwiching \cite{Serna2008} and Hyperboxing \cite{Teichert2013} runs. Both are adaptive Pareto front approximation algorithms that iteratively construct scalarization problems based on previously calculated Pareto optimal solutions. Sandwiching uses a weighted sum approach, whereas Hyperboxing solves Pascoletti-Serafini problems \cite{Pascoletti1984}. The two algorithms are combined to speed up the process since Sandwiching performs more efficient on convex areas of the Pareto front than Hyperboxing. The adaptive iteration runs are repeated until an approximation quality has been reached or the maximum number of solutions has been exceeded. Different algorithms can be applied \cite{Miettinen1999} \cite{Ehrgott2005} \cite{Ruzika2005} \cite{Branke2008}. Their discussion is not the focus of this paper, though. We assume that a good representation of the Pareto front has already been computed and concentrate on the second stage and on the third stage of the process. 

The paper is organized as follows:
In Section \ref{sec:optimization_problem}, we recall the formulation and some notation for the multicriteria optimization problem. Then, in Section \ref{sec:triangulation}, the triangulation technique is introduce{d} and analyzed. Section \ref{sec::navigation} describes the navigation process. In Section \ref{sec::discussion}, the overall approach is discussed in comparison to already established methods. We conclude with Section \ref{sec::summary}.


\section{Multicriteria Optimization Problem} 
\label{sec:optimization_problem}

At the core of this paper is the goal to solve an optimization problem with multiple $n\ge2$ objectives. The general formulation is as follows:
\begin{align*}
\min_{x \in \X} &\ f(x) = \left\{ f_1(x); f_2(x); ... ; f_n(x) \right\} \tag{MOP} 
\end{align*}
The abbreviation MOP stands for multicriteria optimization problem. All objective functions $f_i:\X \to \R$ are minimized with respect to $m \in \N$ continuous parameters. We refer to 
\begin{equation}
\X \subset \R^m
\end{equation}
as the design space of (MOP){, i.e., each dimension of $\X$ represents one parameter.} The design space could be defined by constraints introducing additional functions. Instead, it is simply assumed that these constraints are fulfilled with $\X$ being the set of all feasible parameter values.
{ In this paper, it is assumed that there is no information about the design space other than the initial input for the Pareto front approximation and navigation. The interpolation between feasible solutions could be either feasible or infeasible. So, they are assumed to be feasible.}

\subsection{Optimal solutions}
\label{op_optsol}

The standard concept to define optimality for multicriteria optimization problems is Pareto dominance. Dominance relations between two solutions are induced by ordering cones {\cite{Miettinen1999}}. In this paper, we assume that the ordering cone is the positive orthant. It defines Pareto optimal solutions.

\begin{definition}
Solution {$x^{(1)} \in \X$} is Pareto optimal (efficient or non-dominated) if and only if there exists no other solution {$x^{(2)} \in \X$} such that $f_i(x^{(2)}) \le f_i(x^{(1)})$ for all objectives $i$ and  $f_j(x^{(2)}) < f_j(x^{(1)})$ for at least one objective $j$. {Solution $x^{(1)}$ is dominated otherwise} by any solution that improves in at least one objective. Evaluated Pareto optimal solution $x \in \X$ defines a Pareto point $f(x)$. The set of Pareto points is referred to as the Pareto front.
\end{definition}

The upper bound and the lower bound on the Pareto front indicate the range of efficient designs which can be obtained \cite{Ehrgott2005}.

\begin{definition}
\label{def:ideal_nadir}
The point $y \in \R^n,$ where $n \in \N$ is the number of objective functions, defined by
\[ y_i := \min_{x \in X} f_i(x) \text{ for all } i = 1,...,n \]
is called the ideal point of the multicriteria optimization problem $\min_{x \in X} f(x)$. The point $y \in \R^n$ defined by
\[ y_i := \max_{x \in X \text{ is Pareto optimal}} f_i(x) \text{ for all } i = 1,...,n \]
is called the nadir point.
\end{definition}

It is the goal of multicriteria optimization to compute the Pareto front or at least to find a good representation of the Pareto front. The navigation, i.e. the actual decision making, defines a separate problem.


\section{Pareto front triangulation}
\label{sec:triangulation}

{Assuming} that a good representation of the Pareto front has already been computed, {e.g., the Pareto front is squeezed between some inner and outer approximation with respect to some distance \cite{Teichert2013, Andersson2018} and the distance is small}, let
\begin{equation}
P = \left\{ f(x_1), f(x_2),..., f(x_K) \right\}, K\in\N,
\end{equation}
be the set of Pareto points associated with $K$ Pareto optimal solutions $x_1,x_2,...,x_K$. {The Pareto front could be connected or not. This input is insufficient to infer connectivity. So, the Pareto front is assumed to be connected.} The second stage of our method is to connect the set $P$ in such a way that all possibly feasible {corresponding} solutions in between are represented, as well. These representations are not necessarily Pareto points of the original problem (MOP), but, they can be assumed to be close to the actual Pareto front depending on the approximation error of the Pareto front representation {\cite{Teichert2013, Andersson2018}}. This is specifically useful {for instance} for the optimization of industrial processes represented by graybox models \cite{Asprion2017}.

For the triangulation in the second stage, we assume that the Pareto front is a subset of a $(n-1)$-dimensional manifold and, hence, can be represented by $(n - 1)$-dimensional triangles (simplices). That can be assumed due to the definition of Pareto optimality, which allows us to extend the Pareto front to a surface that can be projected by an homeomorphism along the vector of all ones. The Pareto front triangulation is defined as
\begin{equation}
D = \cup_{i = 1}^N D_i \subset \text{conv}\{P\},
\end{equation}
where $N\in \N$ is the number of triangles/simplices $D_i$. {The approximation $D$ is part of the convex hull spanned by the Pareto points. The triangles/simplices $D_i$ are convex sets that represent the approximating area spanned in between neighboring Pareto points. In non-trivial cases, the triangulation is not unique. This can be seen in Figure 1, where four points are connected by two triangles in two different ways.}

The following Algorithm \ref{tri_alg} computes the Pareto front triangulation based on two ideas. The first idea is {based on the motivating conjecture} that the Pareto front is a surface of dimension $n-1$ which does not fold in the direction of improvement. The second idea is that Pareto front interpolations should be influenced by {the closest neighboring Pareto solutions}. The algorithm consists of three steps illustrated in Figure \ref{fig:triangulation}. {In the first step, fifteen Pareto points $(n = 3)$ are projected in direction $e = (1,1,1)^T$ indicated by the three arrows onto an $e$-plane, a two-dimensional space orthogonal to $e$. In the second step, Delaunay Triangulation algorithm for space dimension $(n-1) = 2$ is performed inside the $e$-plane to compute the triangulation for fifteen projected points. In the third step, the computed triangulation is transferred to the corresponding set of Pareto points (the correspondence is indicated by three arrows) generating an approximation $D$ that consists of sixteen triangles/simplices $D_i$.}

\begin{algorithm}[h]
 \KwData{ Set of Pareto points $P$}
 \KwResult{ Delaunay triangulation $D$ }
 Step 1: With $e$ being the vector of all ones, project $P$ onto an $e$-plane\;
 Step 2: Apply $DelaunayTriangulation$ from \cite{Berg2008} p. 200 to the projected points\;
 Step 3: Apply the resulting triangulation to $P$ defining $D$\;
 \Return D\;
 \caption{Triangulation algorithm}
 \label{tri_alg}
\end{algorithm}

\begin{figure}
\centering
\includegraphics[scale=.3]{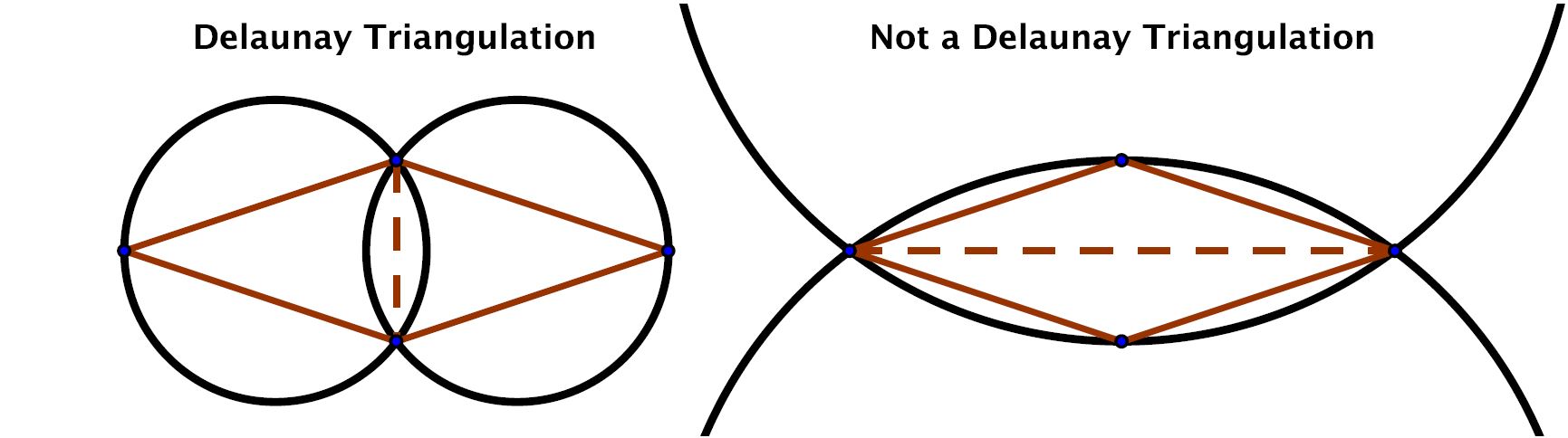}
\caption{Edge legalization in the Delaunay triangulation. {The dashed edge can be flipped to switch from the triangulation in red to the left to the triangulation in red to the right and vice versa.}}
\label{fig:delaunay_egde_flipping}
\end{figure}

\begin{figure}
\centering
\begin{subfigure}{\textwidth}
\centering
\includegraphics[scale=.33]{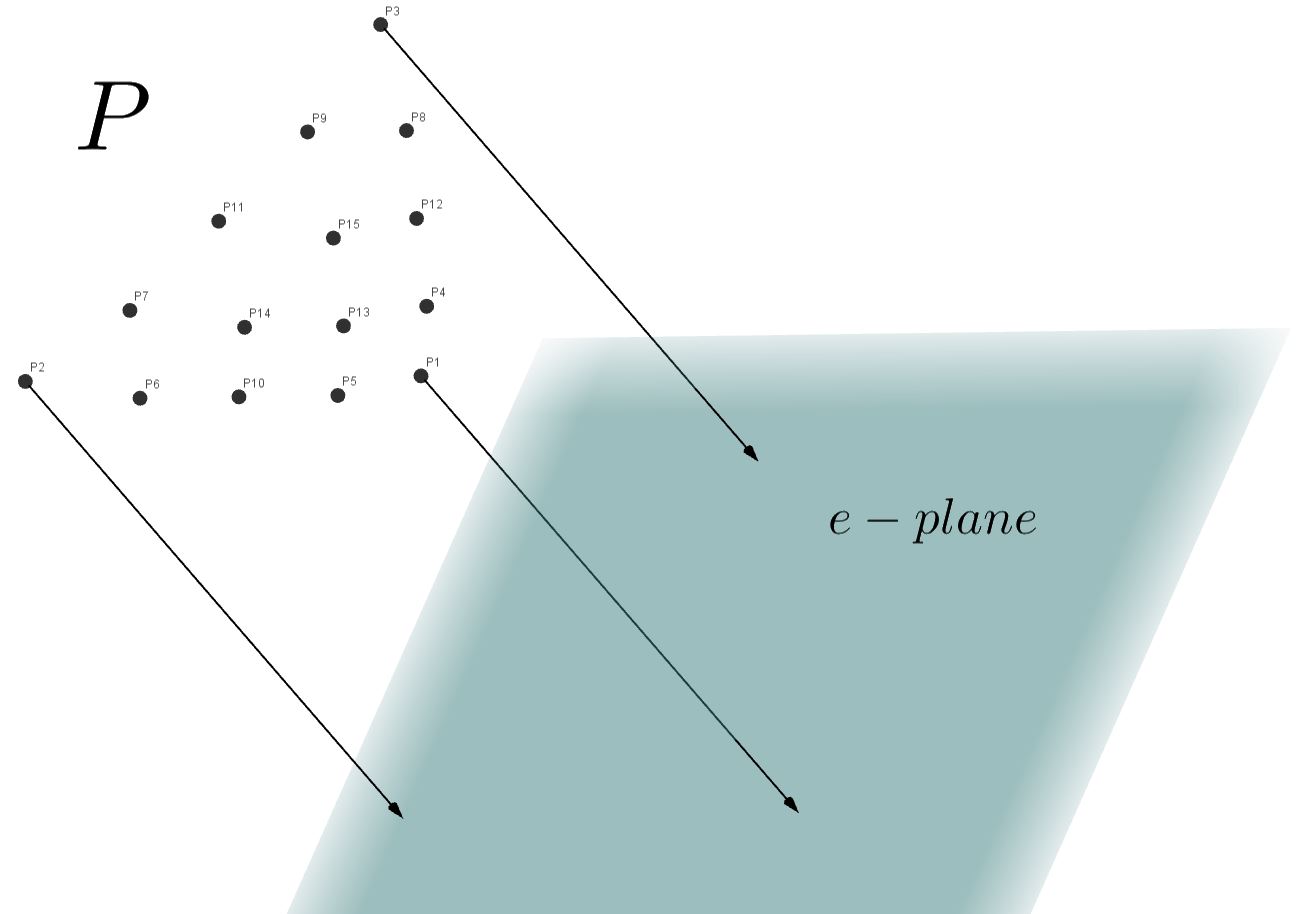}
\caption{Step 1: Projection}
\label{fig:SW_step_1}
\end{subfigure}
\begin{subfigure}{\textwidth}
\centering
\includegraphics[scale=.33]{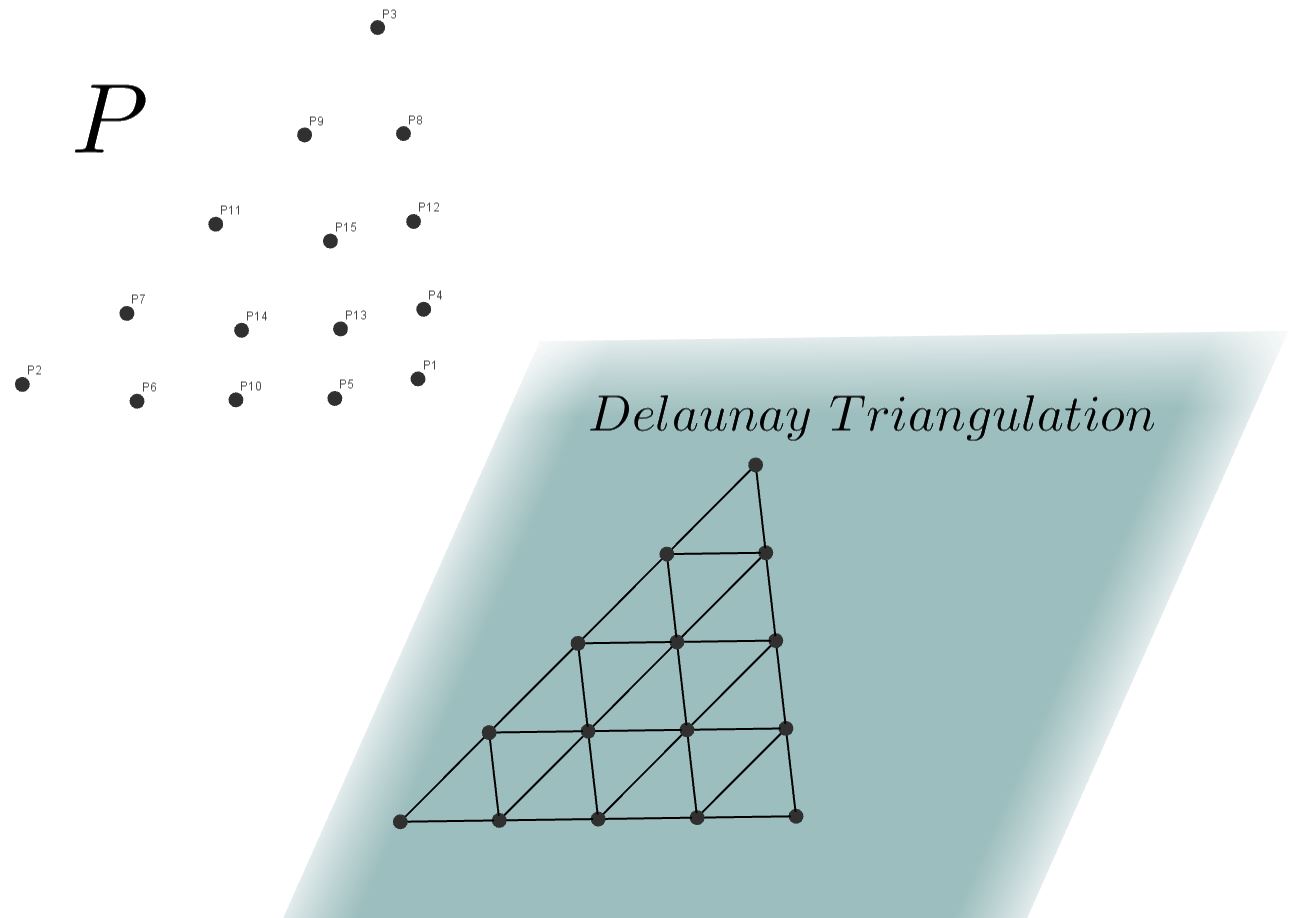}
\caption{Step 2: Delaunay Triangulation}
\label{fig:SW_step_2}
\end{subfigure}
\begin{subfigure}{\textwidth}
\centering
\includegraphics[scale=.33]{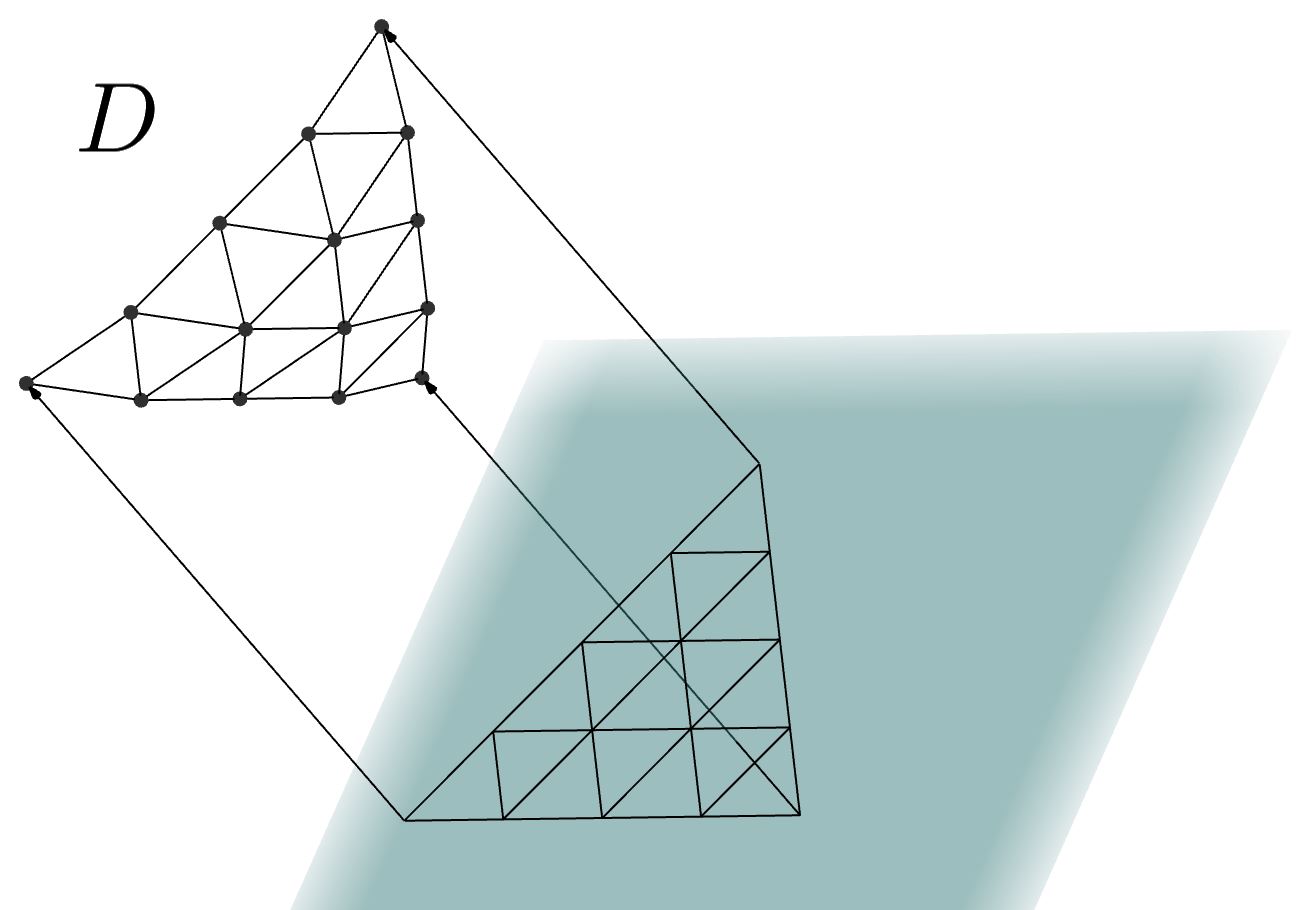}
\caption{Step 3: Transition to triangulation $D$}
\label{fig:SW_step_3}
\end{subfigure}
\caption{Triangulation steps from Algorithm \ref{tri_alg}}
\label{fig:triangulation}
\end{figure}

The purpose of the Delaunay triangulation in the second step of Algorithm \ref{tri_alg} is that approximated areas $D_i, 1 \le i \le N$, are influenced by the closest Pareto points. As depicted in Figure \ref{fig:delaunay_egde_flipping}, a Delaunay triangulation can be constructed by consecutive edge legalizations {\cite{Berg2008}}. In case that an encircling area that is defined by a triangle contains {projected} Pareto points that do not define the triangle, the edge can be flipped {(see, e.g., Figure \ref{fig:delaunay_egde_flipping})}. One of the fastest algorithms to compute a Delaunay triangulation is Quickhull \cite{Barber1996}. It computes the convex hull of the same number of uplifted points in a higher dimension with worst case complexity $O(K\log K)$ for $n \le 3$ and $O\left(K^{\lfloor\frac{n}{2}\rfloor}\right)$ for $n \ge 4$.

\section{Pareto front Navigation}
\label{sec::navigation}

After the Pareto front triangulation has been computed, the third stage of our method involves the decision maker to explore the approximated Pareto front. This navigation process allows two actions, selection and restriction. Selection means that the decision maker selects an objective or parameter to attain certain value. Restriction, on the other hand, {allows filtering the solution space}, i.e. unwanted or uninteresting area of the solution space is removed from consideration.

\begin{figure}
\centering
\begin{subfigure}{\textwidth}
\centering
\includegraphics[scale=.3]{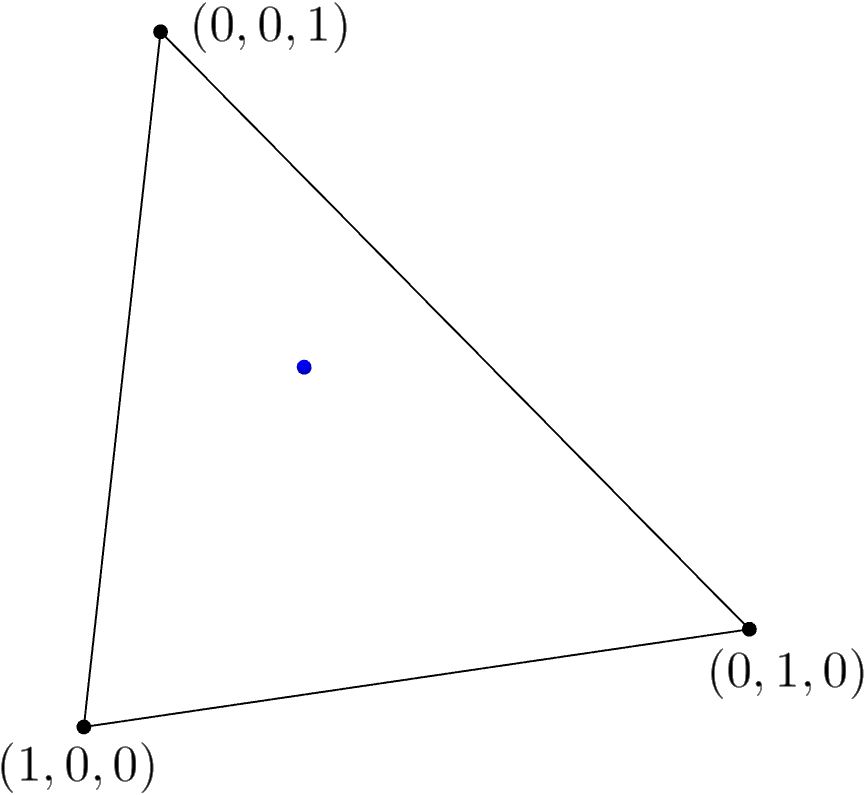}
\label{fig:SW_step_1}
\centering
\includegraphics[scale=.3]{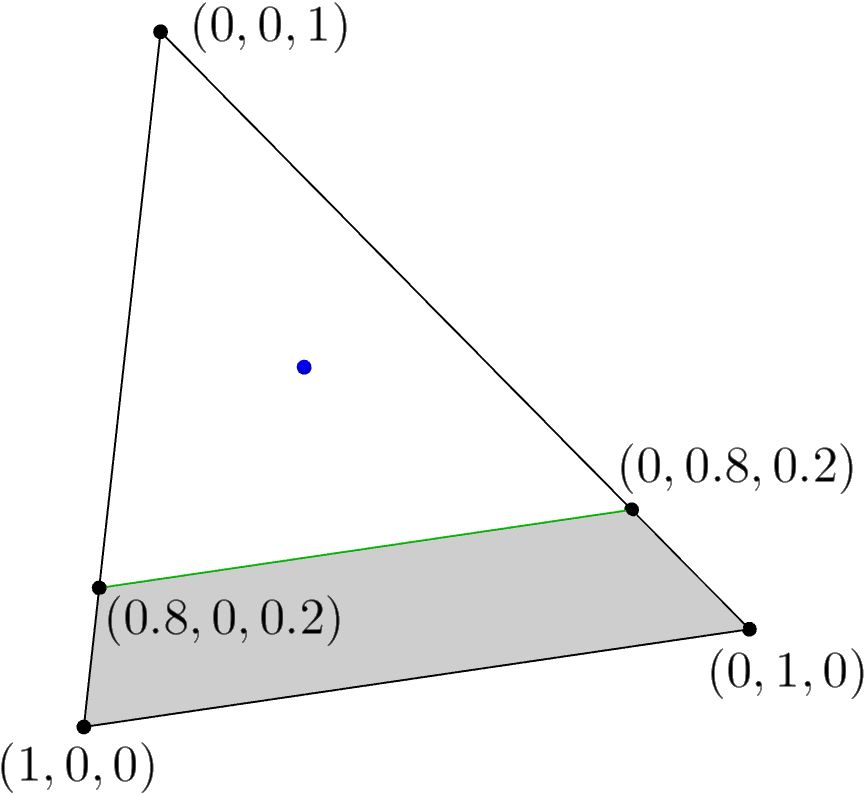}
\label{fig:SW_step_2}
\end{subfigure}

\begin{subfigure}{\textwidth}
\centering
\includegraphics[scale=.63]{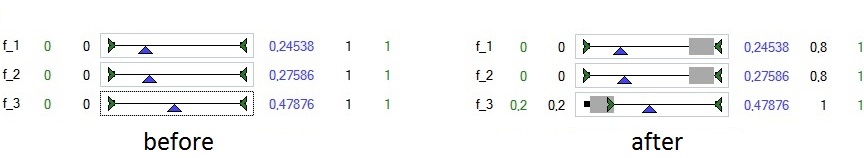}
\label{fig:SW_step_3}
\end{subfigure}
\caption{Navigation on a triangle. Blue color represents the current selection. After active restriction of the solution space in green, the gray area depicts infeasible/unreachable area. {The selectors displayed by the blue arrows can be moved freely inside the implicit/passive bounds. In the beginning, implicit/passive bounds equal the active bounds represented by the restrictors (green arrows). After an active bound is moved, implicit/passive bounds can become tighter extending the infeasible (gray) area on other sliders.} }
\label{fig:restriction}
\end{figure}

Based on the mathematical analysis in \cite{Monz2006}, the selection and restriction mechanism developed at Fraunhofer ITWM was originally intended for intensity modulated radiotherapy (IMRT) planning. A few years later it was implemented for the first time in a software called MIRA - Multicriteria Interactive Radiotherapy Assistant \cite{Monz2008}. Back then the navigation has been visualized by a navigation star centering several sliders, each representing a different objective. Later on, the representation has evolved into a cleaner visualization where each slider is placed one below the other representing either {an} objective or {a} parameter. Figure \ref{fig:restriction} describes a navigation tool that is implemented in many applications at Fraunhofer ITWM.

For illustrative purposes, let us assume that each slider represents an objective. In Figure \ref{fig:restriction}, the decision maker is allowed to change the selected values of three objectives by dragging the blue triangles. They are referred to as selectors. The decision maker is also allowed to actively reduce the range for all three objectives. These so called restrictors or active bounds are indicated by the green arrows. In this example, we start with a {Pareto front approximation} spanned between the points $(1,0,0)^T$, $(0,1,0)^T$ and $(0,0,1)^T$. The selected solution lies inside a triangle. Then, the decision maker moves the left restrictor of the third objective and implicitly reduces the feasible area of the solution space. As a consequence, the so called implicit or passive upper bounds of the first two objectives reduce from $1$ to $0.8$. {During the navigation, selection and restriction are often combined until a desirable Pareto point/solution is reached. Recall, that the result is an approximation. The accuracy of the result depends on the approximation quality that was used to generate the Pareto set. In order to assure accuracy or feasibility, it is possible to use the result as a warm start for an additional post-processing step that finds the closest feasible Pareto point to the result, e.g., by solving a Pascoletti-Serafini problem \cite{Pascoletti1984}.}

In comparison to the interactive solution process of NIMBUS \cite{NIMBUS}, the classification into five different objective classes could be translated by the sliders as follow{s}:
\begin{itemize}
\item $I^<$: An objective is to be improved (by dragging the selector of a slider to the left). The active upper bound is set to the current selector value. {E.g., in order to improve the objective $f_1$ on the simplex in Figure \ref{fig:restriction}, the user can drag the corresponding blue selector to the left up to the value $0$. The improvement of the objective $f_1$, however, deteriorates the remaining objectives $f_2$ and $f_3$ according to the condition $f_1 + f_2 + f_3 = 1$.}
\item $I^\le$: An objective is improved by dragging the selector of a slider to the desired aspiration level. The active upper bound is set to the aspiration level. {E.g., after the objective $f_1$ in Figure \ref{fig:restriction} is set to the value $0.1$, the active upper bound represented by the green arrow to the right of the slider is also set to the same value (by drag and drop). This improvement, though, comes with the cost that $f_2 + f_3 \ge 0.9$.}
\item $I^=$: The objective value must not deteriorate. The active upper bound is set to {the} current selector value. {(see Example for $I^\le$)}
\item $I^\ge$: The objective value is allowed to deteriorate. The active upper bound is set to the allowed worst case value. {E.g,. setting the upper bound of objective $f_1$ to the worst case value 0.5 in Figure \ref{fig:restriction} implies the condition $f_2 + f_3 \ge 0.5$, i.e., we allow our desired value 0.1 to deteriorate to the value 0.5.}
\item $I^{<>}$: The objective is allowed to change freely. The restrictors of the slider stay untouched. {I.e., if the active bounds of an objective represented by the green arrows in Figure \ref{fig:restriction} are not set, the blue sliders can be moved freely as long as they do not cross the implicitely restricted area denoted by the gray boxes inside the slider.}
\end{itemize}
If the restriction imposed by the upper bounds can not be set with the sliders, i.e., the solution space with these restrictions is empty, then, the user has to modify (deteriorate) the aspiration levels until a feasible solution can be selected. The objectives are improved by dragging the selector to the left. In order to fix the improvement in any objective, the active upper bound is moved to that value and the user continues with the next objective to be improved. If any objective can not be improved more than indicated by the passive bounds of the sliders, the user can adjust some active bounds and observe by the change of the passive bounds how much improvement can be gained by these adjustments.

In the following subsections, we formulate the problems that are solved during the navigation process each time a selector or restrictor is moved. We distinguish between the selection and restriction of objectives and parameters since the Pareto front is computed in the objective space and not in the design space. For this reason, a variety of different cases has to be taken into consideration. For simplification reasons, we no longer distinguish between the terms triangulation, Pareto front and Pareto front approximation.


\subsection{Selection of an objective value}
\label{subsec:objective_selection}

\begin{figure}
\centering
\includegraphics[scale=.4]{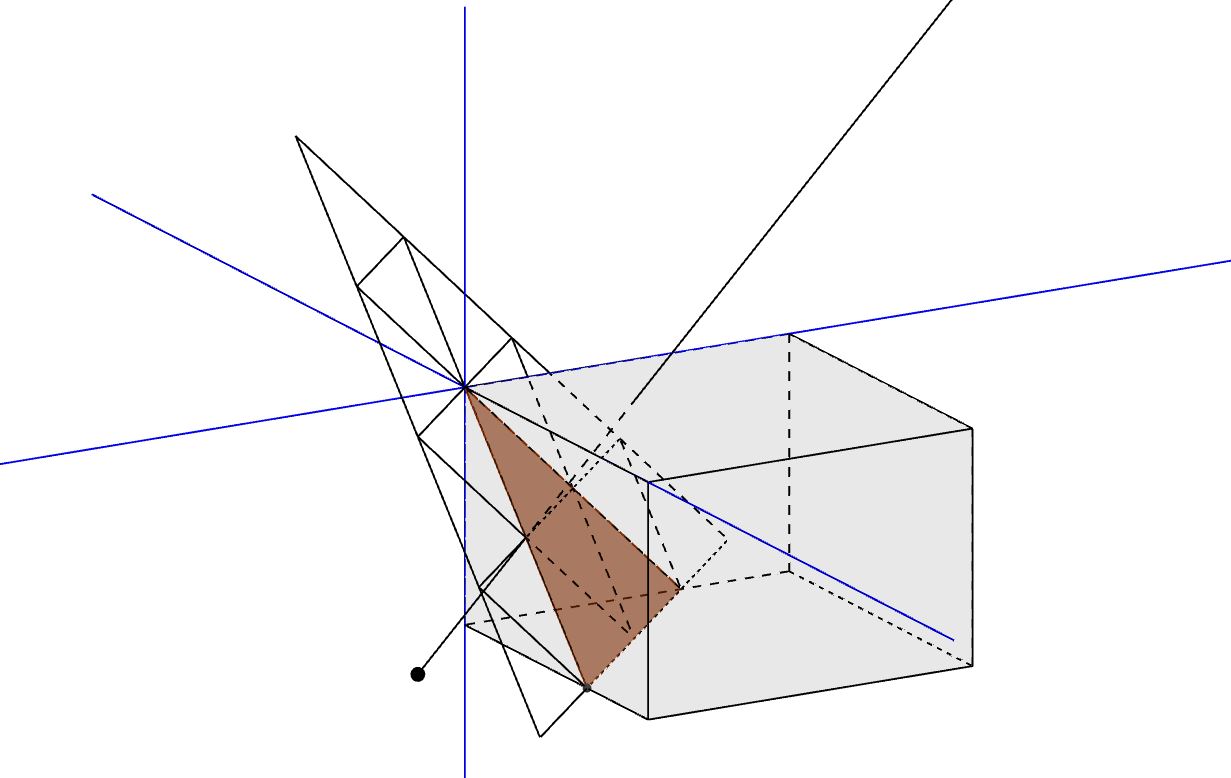}
\caption{Ray tracing on a triangulation using a kd-tree data structure. {The triangulation consists of 16 simplices. Three splitting hyperplanes split the area into 8 boxes. The three blue lines indicate the intersections between each of the two splitting hyperplanes. The highlighted area illustrates the interior of one box containing 4 red simplices. Ray tracing shoots a ray starting at a reference point displayed by the black dot. This ray hits the highlighted area and, thus, the box containing the four triangles. It remains to be checked whether this ray also intersects with one of the red simplices. Simplices inside boxes, that are not hit by the ray, are not checked for intersection.}}
\label{fig:kd_tree_raytracing}
\end{figure}

The selection is a grab-move-release mechanism that allows the decision maker to explore the feasible area of the Pareto front approximation. {Desirably, the interaction should be fast and reactive not only} after the release action of a selector but also during the {movement}. In order to guarantee a good user experience, the selection mechanism was implemented based on a ray tracing technique from computer graphics \cite{Hapala2011} that makes use of a kd-tree data structure to quickly find triangles representing the Pareto front. The main idea {of} this approach is illustrated in Figure \ref{fig:kd_tree_raytracing}.

During the {movement} of a selector, the {corresponding value} of the selected point changes and creates a reference point that {can be} farther away from the triangulation. In order to get back onto the triangulation, a ray is shot into the direction $e_{-k}$, a vector of all ones but zero in the $k$-th entry, where $k\in\{1,...,n\}$ is the modified selector/objective. The reference point is projected back onto the Pareto front. By definition of Pareto optimality, we move to a new solution with the same objective value as chosen by the selector and the most favorable distance to the reference point with respect to the other objectives. Notice that for convex Pareto front approximations this approach is analogous to the Pascoletti and Serafini approach in \cite{Monz2008}. {For the non-convex case, the approximations are not guaranteed to be non-dominated which is the main difference to the Pascoletti and Serafini approach.}

The performance of ray tracing is improved by using a kd-tree data structure that stores the triangulation $D$ defined in Section \ref{sec:triangulation}. Given the ideal point and the nadir point of the triangulation {(see Definition \ref{def:ideal_nadir})}, the corresponding box is split in each dimension (possibly in two equal parts) to create $2^n$ nodes/boxes in the next level of the kd-tree. In Figure \ref{fig:kd_tree_raytracing}, the split creates $8$ boxes. The highlighted area represents {the interior of} one box containing four red triangles. This procedure can be recursively repeated until each node/box of the kd-tree covers/stores a small number of triangles to check whether the ray intersects the triangle. So, it suffices to check triangles inside {all} boxes that are intersected by the ray. Thus, the selection time complexity for ray tracing can be reduced to $O(\log(N))$. 

{The intersection of the ray with a box can be computed in linear time with respect to $n$. Let ray $r$ be defined by 
\[ r(\mu) = p + \mu e_{-k}, \mu \ge 0, \] 
where $p$ is the starting point and $e_{-k}$ the direction of the ray. Then, by determining in each dimension the range of $\mu$ for which the box and the ray intersect, the intersection of all ranges determines that the box is either missed in case of emptiness or hit otherwise. E.g., for two objectives, let us assume that the user starts at the approximated Pareto point $(2,2)^T$. Then, the selector of the first objective is moved to value $1$, i.e., $p = (1,2)^T$ and 
\[ r(\mu) = (1,2)^T + \mu (0,1)^T, \mu \ge 0. \] 
The box $B = [0,2]\times[2,4]$ can be tested for intersection by checking the intersection for each dimension. With respect to the first dimension, the ray intersects for all $\mu \ge 0$, since $r_1(\mu)=1\in[0,2]$. With respect to the second dimension, $r_2(\mu) \in [2,4] \iff \mu\in[0,2]$. I.e., the ray intersects the box B for $\mu\in[0,\infty)\cap[0,2] = [0,2]$.

The intersection of a ray with a simplex is the solution to the linear equation system defined by the ray and the plane extending the $(n-1)$-dimensional triangle. If the result is a convex combination of the supporting Pareto points, then the ray should intersect the interior of the triangle. Going back to the two-dimensional example, assume that the box $B$ contains a simplex defined by the Pareto points $q^{(0)} = (2,2)^T$ and $q^{(1)} = (0,3)^T$. Then, the resulting linear equation to be solved is:
\begin{align*}
1 + 0\cdot\mu =& 2 + (0-2)\cdot\lambda_1 \\
2 + 1\cdot\mu =& 2 + (3-2)\cdot\lambda_1
\end{align*}
The solution to this equation is $\mu = 0.5$ and $\lambda_1 = 0.5$. Since $\lambda_1 = 0.5 \in [0,1]$, the simplex is intersected in the interior of the simplex spanned by the Pareto points $q^{(0)}$ and $q^{(1)}$.
}

\begin{figure}
\centering
\includegraphics[scale=.4]{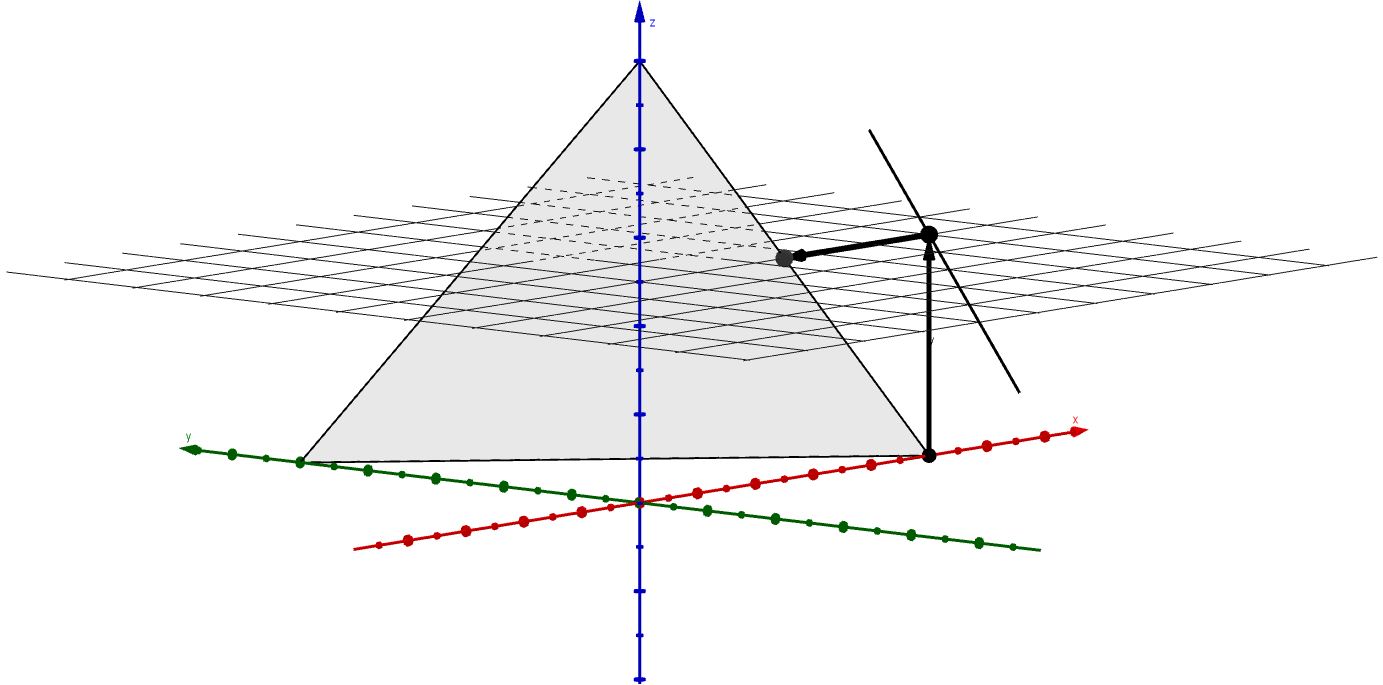}
\caption{Ray tracing misses the triangulation. {The user starts in the lower right corner of the triangle. The first black arrow indicates the change of the blue objective value ending in a reference point (black dot outside the triangle). At the reference point, two rays are shot in both directions inside the meshed area maintaining the blue objectives selected value and miss. So, the reference point is projected to the closest feasible approximated Pareto point with the chosen selector value. The user is navigated according to the second black arrow to this point.}}
\label{fig:raytracing_misses}
\end{figure}

Sometimes ray tracing may not hit the triangulation. This can happen for different reasons. One reason may be that the triangulation is in the opposite direction. Then, we simply use the opposite direction, i.e, the vector $-e_{-k}$. Another reason may be that the ray is next to the border of the triangle (see Figure \ref{fig:raytracing_misses}). In that case, a different auxiliary problem has to be solved.

When ray tracing misses {all} triangle{s}, the closest solution to the reference point does not necessarily have the same value as the moved selector. As depicted in Figure \ref{fig:raytracing_misses}, the distance has to be computed on a plane/set of points with the same selector values. Given the reference point $p \in \R^n$ and the triangle $\Delta q^{(0)} q^{(1)}... q^{(n-1)}$, $q^{(0)} \in P$, $q^{(j)} \in P$ for all $j \in \{1,...,n-1\} $, we solve the following problem:
\begin{align*}
\min_d\ & d \\
\text{s.t. } & d \ge |y_i - p_i| \text{ for all } i \neq k \\
& y_k = p_k \\
& y_i = q_i^{(0)} + \sum_j \lambda_j \left( q_i^{(j)} - q_i^{(0)} \right) \text{ for all } i \\
& \sum_j \lambda_j \le 1, \lambda_j \ge 0 \text{ for all } j
\end{align*}
which can be reformulated {as} the linear problem:
\begin{align*}
\min_d\ & d \\
\text{s.t. } & \sum_j \lambda_j \left( q_k^{(j)} - q_k^{(0)} \right) = p_k - q_k^{(0)} \\
 & \sum_j \lambda_j \left( q_i^{(j)} - q_i^{(0)} \right) - d \le p_i - q_i^{(0)} \text{ for all } i \neq k \\
& \sum_j \lambda_j \left( q_i^{(0)} - q_i^{(j)} \right) - d \le q_i^{(0)} - p_i \text{ for all } i \neq k \\
& \sum_j \lambda_j \le 1, \lambda_j \ge 0 \text{ for all } j
\end{align*}

A linear problem can be solved in polynomial time \cite{Jar2004}. Still, in comparison to the previously mentioned operations, this is quite expensive. If, however, the selectors are mostly used to strictly improve the corresponding objectives, ray tracing will hardly miss the triangulation. Also, the kd-tree structure and storing the nadir point and ideal point information for each triangle can be further used to precompute a lower bound on the distance which will hint to the closest triangles and immensely reduce the number of linear programs to be solved.


\subsection{Selection of a parameter value}
\label{subsec:parameter_selection}

The selection of a parameter value is more complex than the selection of an objective value. Ray tracing can not be applied in the objective space because no ray in the objective space can guarantee that the selected parameter value will stay the same once the triangulation is hit. {After transferring a triangulation $D$ from Section \ref{sec:triangulation} to the design space,} ray tracing can not be applied in the design space because it was originally created in the objective space and may behave more chaotic in the design space (by overlapping and folding the triangulation{, i.e., ray tracing may hit several simplices at once}) not to mention the difficulties that may occur when the design space has very different dimension than the objective space.

\begin{figure}
\centering
\includegraphics[scale=.3]{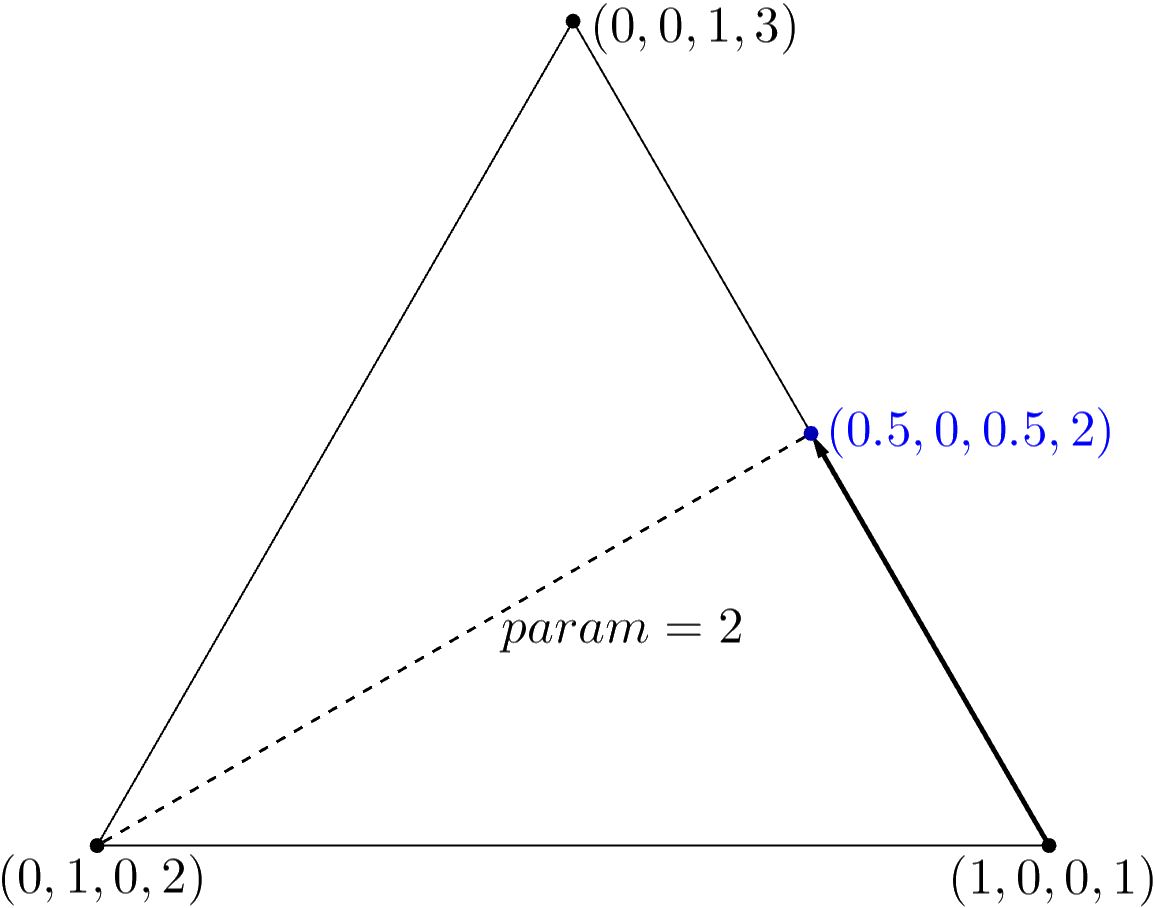}
\caption{Parameter Navigation moving the selector from value 1 to 2: {The user starts with Pareto point $(1,0,0)^T$ and parameter value $1$ at the lower right corner. The interpolated Pareto points with parameter value $2$ are displayed by the dashed line. The closest Pareto point is $(0.5,0,0.5)^T$. The arrow indicates the navigation trajectory for dragging the parameter selector.}}
\label{fig:parameter_navigation}
\end{figure}

The idea is to use an additional kd-tree structure in the design space, i.e. the triangles $D$ are also stored in a different tree specifically generated for quick selection of triangles relevant for parameter selection. So, instead of intersecting a ray with boxes of the kd-tree, we check {the} boxes for selected parameter slider values, which can instantly be determined from the range of a box. {The triangulation $D$ is naturally extended to the parameter space, i.e., the triangulation from the objective space is transferred to the corresponding solutions in the design space $\X$. Given the smallest box containing the triangulation in the design space, the kd-tree structure is constructed analogously to the objective space by splitting the space into smaller boxes at each parameter dimension (see Section \ref{subsec:objective_selection}). Instead of ray tracing, the boxes and simpices are checked whether the navigated parameter value can be achieved by the box or triangle. E.g., if the second parameter selector is moved to some value $v \in \R$, the box $C = [0,1]\times[2,3]$ in the design space contains the navigated parameter value if and only if $v \in [2,3]$. Similarly, simplices are checked in terms of their enclosing bounding boxes (analogous to the nadir point and ideal point information of a simplex in the objective space). After this a priori selection of relevant triangles, each triangle is tested for the closest approximated Pareto point/solution with selected parameter value. The user is navigated to the closest result.

In this paragraph, the technical computation of the closest approximated Pareto point/solution is explained. }Let us assume that the objective space information is naturally extended by the corresponding parameter information. Given the reference point $p \in \R^{n + m}$ and the triangle $\Delta q^{(0)} q^{(1)}... q^{(n-1)}$, $q^{(0)} \in P\times\X$, $q^{(j)} \in P\times\X$ for all $j \in \{1,...,n-1\} $, we select the closest solution in the approximated Pareto set with selected value of {parameter $k \in \{n+1,...,n+m\}$} by solving the following linear program:
\begin{align*}
\min_d\ & d \\
\text{s.t. } & \sum_j \lambda_j \left( q_i^{(j)} - q_i^{(0)} \right) - d \le p_i - q_i^{(0)} \text{ for all } i  \in \{1,...,n\} \\
 & \sum_j \lambda_j \left( q_i^{(0)} - q_i^{(j)} \right) - d \le q_i^{(0)} - p_i \text{ for all } i  \in \{1,...,n\} \\
& \sum_j \lambda_j \left( q_k^{(j)} - q_k^{(0)} \right) = p_k - q_k^{(0)} \\
& \sum_{j} \lambda_j \le 1, \lambda_j \ge 0 \text{ for all } j
\end{align*}
Consider the example in Figure \ref{fig:parameter_navigation}. Starting with a triangle spanned by the Pareto points $(1,0,0)^T, (0,1,0)^T$ and $(0,0,1)^T$ with parameter slider values $1$, $2$ and $3$, respectively (other parameters are omitted here), the selected value $1$ is increased by $1$. Since the selected point is not moved in the objective space, we have to find the next proximate point on the simplex with parameter value 2. Possible candidates are represented by the {dashed} segment line. Point $(0.5, 0, 0.5)^T$
with parameter value 2 is closest to the lower right corner. Therefore, point $(0.5, 0, 0.5)^T$ is navigated to next.


\subsection{Restriction}

Whenever a restriction is set {to an} objective function or {to a} parameter values, some of the attainable ranges for navigation may be restricted as well (see for example Figure \ref{fig:restriction}). These attainable passive ranges have to be separately recalculated each time the restrictors are moved, so that the user understands how restriction changes the range of the feasible solution space. Unfortunately, the computation of passive bounds, {i.e., bounds that can be implicitely inferred from active bounds (see, e.g., gray area in Figure \ref{fig:restriction})}, is costly. However, once again, the kd-tree structures in the objective space and in the design space can be reused to quickly find triangles that are influenced by the corresponding restrictor {movement}.

\begin{figure}
\centering
\begin{subfigure}{\textwidth}
\centering
\includegraphics[scale=.25]{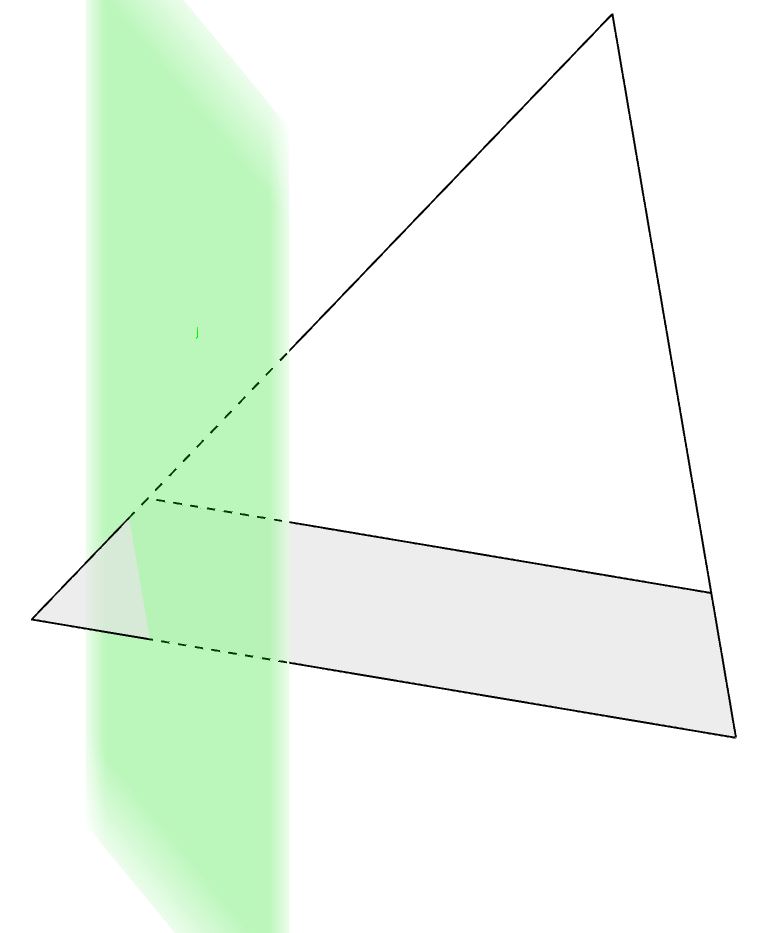}
\centering
\includegraphics[scale=.25]{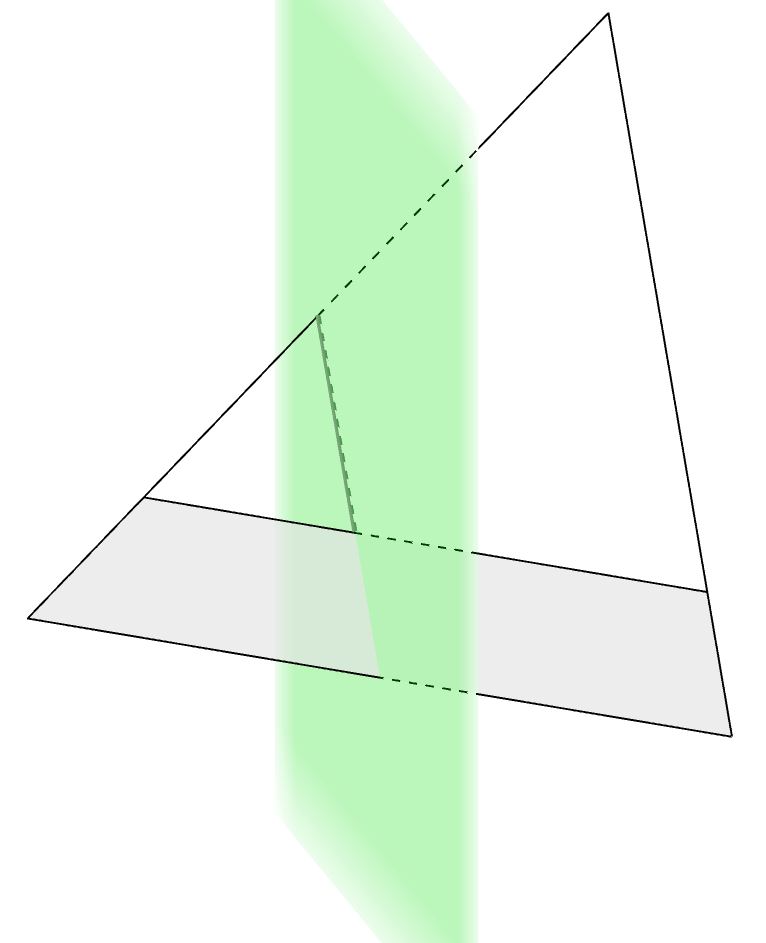}
\end{subfigure}
\caption{{Restriction of an already restricted triangle/simplex with the restricted area displayed in gray and the unrestricted interior area displayed in white. The hyperplane in green further restricts the area to the left of the hyperplane. The triangle to the right is influenced. The triangle to the left, however, is not influenced. }}
\label{fig:restriction_of_restricted_area}
\end{figure}

After restrictors have been moved {(green arrows in Figure \ref{fig:restriction})}, some triangles are intersected by the restricting hyperplane (see Figure \ref{fig:restriction_of_restricted_area}). {The restricting hyperplane is either the plane defined by the moved objective restrictor, e.g., if the third objective's (out of $n = 5$ objectives) lower restrictor has been moved to value $v\in \R$, the restricting hyperplane in the objective space is $\{x\in\R^5\mid x_3 = v\}$, or (analogously) the restricting hyperplane is defined in the design space by the moved parameter restrictor. The intersected simplices redefine the feasible area. The feasible part of the approximated Pareto front and the implicitly infeasible area (gray area inside the sliders in Figure \ref{fig:restriction}) can be entirely computed from the simplices influenced by the restricting hyperplane. For this reason, the feasible area inside these simplices has to be reevaluated. Influenced simplices} are squeezed inside a bounding box $[a,b]$ represented by the {(ideal)} point $a\in\R^{n + m}$ and the {(nadir)} point $b\in\R^{n + m}$ of the feasible interior area of the triangle. This bounding box information is crucial for the evaluation of {implicit/passive bounds} and should be stored for later restrictor move actions to check whether the bounding box has to be updated. Given the active lower bounds and upper bounds $l,u\in\R^{n+m}$ and the simplex $\Delta q^{(0)} q^{(1)}... q^{(n-1)}$, $q^{(0)} \in P\times\X$, $q^{(j)} \in P\times\X$ for all $j \in \{1,...,n-1\} $, we have to compute each entry of the bounding box by solving the following linear programs:
\begin{align*}
\min_{a_k}\ & a_k \\
\text{s.t. } & a_k = q_k^{(0)} + \sum_j \lambda_j \left( q_k^{(j)} - q_k^{(0)} \right) \\
& \sum_j \lambda_j \le 1, \lambda_j \ge 0 \text{ for all } j \\
& a \in [l,u]
\end{align*}
which is equivalent to:
\begin{align*}
\min_{a_k}\ & a_k \\
\text{s.t. } & \sum_j \lambda_j \left( q_k^{(j)} - q_k^{(0)} \right) - a_k = - q_k^{(0)} \\
& l_i - q_i^{(0)} \le \sum_j \lambda_j \left( q_i^{(j)} - q_i^{(0)} \right) \le u_i - q_i^{(0)} \text{ for all } i \\
& \sum_j \lambda_j \le 1, \lambda_j \ge 0 \text{ for all } j
\end{align*}
for the lower bound entries and an analogous $\max_{b_k}\ b_k$ formulation for the upper bound entries for all $1 \le k \le n+m$.
{Every time the active bounds represented by the slider restrictors in Figure \ref{fig:restriction} are moved, the passive bounds are reevaluated from the bounding box information of influenced triangles and adjusted accordingly to be represented in the slider in terms of gray boxes. Notice that passive/implicit bounds are tighter that the active bounds.}

\subsection{Selection on restricted triangulation}

The restriction of Pareto front influences the selection process. The selected solutions are not allowed to lie inside the infeasible area of the triangulation. 

The ray tracing mechanism described in Subsection \ref{subsec:objective_selection} has to be adapted. After intersection with {the} triangulation $D$, the solution is checked for feasibility. If active restrictions of all sliders are fulfilled, the solution is valid. Otherwise, an auxiliary problem similar to the linear program from Subsection \ref{subsec:objective_selection} has to be solved. The only change is an extension by additional linear constraints
\[ l_i - q_i^{(0)} \le \sum_j \lambda_j \left( q_i^{(j)} - q_i^{(0)} \right) \le u_i - q_i^{(0)} \text{ for all } i, \]
where $l,u\in\R^{n+m}$ are the active lower bounds and upper bounds set by the decision maker.
The parameter selection mechanism in Subsection \ref{subsec:parameter_selection} does not make use of ray tracing. It is solely adapted by adding the same constraints to the corresponding auxiliary linear programs to be solved.

\begin{figure}[!htb]
\centering
\includegraphics[scale=.3]{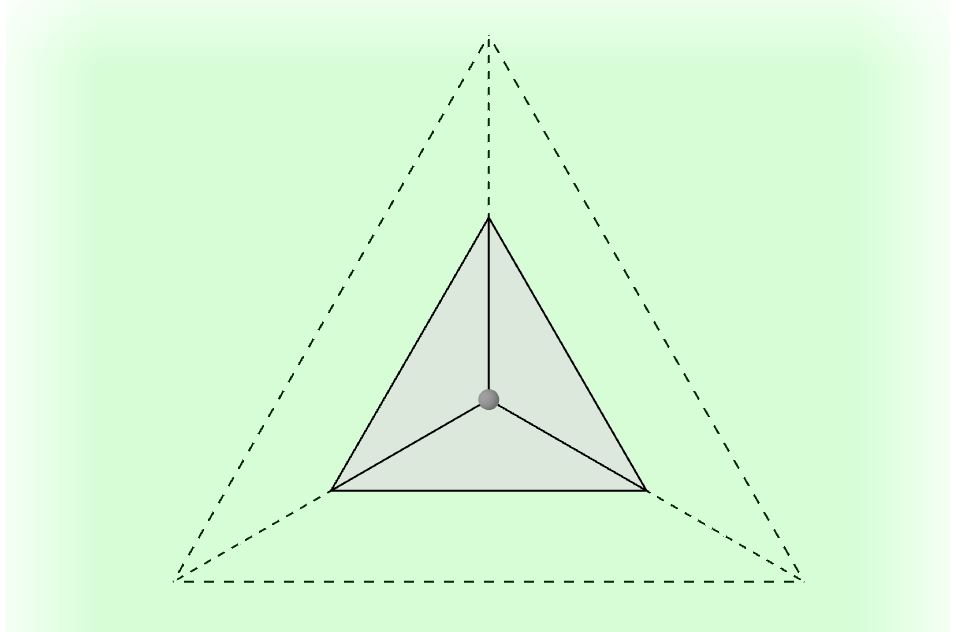}
\caption{Parameter restriction {may generate} holes in the triangulation: {The center Pareto point has a low parameter value. The three Pareto points in the corners of the triangulation have the same high parameter value. The three triangles (n = 3) in the objective space are restricted by the parameter's lower restrictor move which creates a hole in the Pareto front triangulation. The resulting infeasible area in the center is displayed by gray color. }}
\label{fig:parameter_holes}
\end{figure}

Parameter restriction may cause a different issue. As depicted in Figure \ref{fig:parameter_holes}, the center area may be cut off by the restricting parameter hyperplanes. It is possible to create examples with almost as many holes in the triangulation as the number of {initial} Pareto points. Then, some areas on the sliders may not be selected. Other times, the selector may not move continuously as the user experiences abrupt selector jumps. For the user of the navigation tool described in Figure \ref{fig:restriction} this creates a bad experience. To avoid this, we decided to check feasibility of ray tracing intersections solely with respect to the restriction for objectives. If parameter restrictions are not fulfilled, the corresponding entries are modified to the next feasible values and the selector is turned {from blue to} gray to indicate that the selected point is in fact infeasible (see Figure \ref{fig:parameter_holes}). {E.g., assume that in Fugure \ref{fig:parameter_holes} the parameter values range from 0 to 1. Then, setting the parameter's lower restrictor to 0.5 generates a hole in the triangulation. When ray tracing is used and the ray intersects with the gray area in Figure \ref{fig:parameter_holes}, the objective values are still valid with respect to objective slider restrictors. The corresponding parameter value of the approximated Pareto point, which is below the value 0.5 and, thus, infeasible, is adjusted to the value 0.5 and the selectors are turned from blue to gray to indicate this infeasibility. This way, the navigation performs smoother (without jumps) when holes are involved. }

\section{Discussion}
\label{sec::discussion}

Our method can be {best} compared to an already established interactive multicriteria optimization method PAINT \cite{Hartikainen2012} which was specifically developed for nonconvex problems. It is often related with NIMBUS. PAINT (PAreto front INTerpolation) is an approximation method that is used in the second stage. The main advantage of PAINT is that a set $P$ of Pareto optimal solutions computed in the first stage is connected in such a way, that all interpolated solutions in between stay Pareto optimal. This is achieved by triangulating the set $P$ and, then, constructing a sub-complex by cutting off simplices which either dominate or get dominated by any solution in $P$. The resulting subset of the convex hull of the set $P$ is called inherently nondominated \cite{Hartikainen2011}. Due to the comparison of simplices, PAINT is very time consuming at higher dimensions. For the DTLZ2 test problem \cite{DLTZ} with 20 randomly generated points in $\R^4$, PAINT takes 111.2 seconds. In comparison, our second stage generates an approximation of 100 randomly generated points in $\R^{10}$ in about the same time. {Table \ref{tab:time_complexity} displays the time consumption arising from Delaunay triangulation. For more than ten objectives, this triangulation becomes quite slow. So, it is not recommended to apply this method for more than 10 objective functions. Fortunately, the number of parameters does not matter since the triangulation from the objective space is transferred to the design space.  }

Another feature of PAINT is how holes in the Pareto front are dealt with (see an application to the three-objective Viennet's test problem \cite{VlENNET}). PAINT detects and takes account of holes, such that during navigation (stage three), sudden jumps can occur. On the one hand, this is an advantage since, from a theoretical point of view, a Pareto front is not necessarily connected. On the other hand, this is a disadvantage since, from a decision maker point of view, an interactive exploration of an approximation with holes will involve rapid jumps in the objective space. Therefore, our navigation approach deals with holes in a different way. 

\begin{table}[htbp]
  \begin{center}
  \caption{{Time consumption of the triangulation step for $n$ objectives (in seconds).}}
    {\begin{tabular}{lcccc}
      \toprule
      number of Pareto points & 10 & 100 & 1000& 10000 \\
      \midrule 
      $n = 6$ & $\sim0.05$ & $\sim0.33$ & $\sim7.72$ & $\sim133.42$ \\
      $n = 7$ & $\sim0.06$ & $\sim1.83$ & $\sim60.88$ & $>1000$ \\
      $n = 8$ & $\sim0.05$ & $\sim8.49$ & $\sim557$ & $>1000$ \\
      $n = 9$ & $\sim0.05$ & $\sim37.83$ & $>1000$ & $>1000$ \\
      $n = 10$ & $\sim0.05$ & $\sim183.77$ & $>1000$ & $>1000$ \\
      \bottomrule
    \end{tabular} }
    \label{tab:time_complexity}
  \end{center}
\end{table}

{The latest successor of PAINT is the method called PAINT-SiCon \cite{Hartikainen2015}. It makes use of additional sample point information in the design space to identify components of the Pareto set (in the design space) that define connected areas of the Pareto front and allow a more accurate approximation of the parameter values during the Pareto front interpolation. However, it does not improve the time complexity of PAINT, which is still used as a sub-algorithm. In this paper, it is assumed that the Pareto set is one component. Also, there is no sample point information beside the non-dominated solutions. The extension to several components (already implemented in the current navigation tool) will be part of future work.}

The navigation process presented in this paper also differs from NIMBUS \cite{NIMBUS}. NIMBUS asks the user what objectives should be improved, at the cost of which objectives and up to what extent \cite{Ham}. The decision maker has to invest more time into the decision process and to reconsider his preferences after every step. Our approach differs by scanning the Pareto front approximation until a satisfactory solution has been found. NIMBUS solely defines an interactive selection process. Our method allows filtering of both the decision space and the objective space and informs the decision maker about implicit changes in restricted dimension ranges. {Tested for no more than 10 objective functions and an arbitrary number of parameters, the slider navigation is smooth in real time.}


\section{Summary and Future Work}
\label{sec::summary}

In this paper, we have introduced a triangulation for Pareto front representations specifically designed to properly navigate on an approximation of the Pareto front with emphasis on a responsive, comprehensive and quick experience for the decision maker. As discussed in Section \ref{sec::discussion}, {the proposed method} has both advantages and disadvantages {with respect to} some existing methods. But, in our opinion, the advantages overweight due to applicability to problems with up to 10 objectives and an arbitrary number of parameters. The proposed method has already been tested by our partners from the chemical industry and found successful.

Our next goal is to improve this method specifically with respect to monotonous objective functions on convex design spaces, which is often the case in industrial applications.

\bibliographystyle{unsrt}  


\begin{thebibliography}{1}

\bibitem{Allmendinger2017} Allmendinger R., Ehrgott M., Gandibleux X., Geiger M. J., Klamroth K., Luque M.: \textit{"Navigation in multiobjective optimization methods"}. "Journal of Multi-Criteria Decision Analysis", 24(1-2):57-70, 2017.

\bibitem{Andersson2018} Andersson T., Nowak D., Johnson T., Mark A., Edelvik F., K\"ufer K. : \textit{"Multiobjective Optimization of a Heat-Sink Design Using the Sandwiching Algorithm and an Immersed Boundary Conjugate Heat Transfer Solver"}. "ASME. J. Heat Transfer", 140(10), 2018.

\bibitem{Asprion2017} Asprion N., Böttcher R., Pack R., Stavrou M.-E., Höller J., Schwientek J., Bortz M.: \textit{"Graybox Models - New Opportunities for the Optimization of Entire Processes"}. "27th European Symposium on Computer Aided Process Engineering", 40(A):97-102, 2017.

\bibitem{Barber1996} Barber C. B., Dobkin D. P., Huhdanpaa H.: \textit{"The Quickhull Algorithm for Convex Hulls"}. "ACM Trans. Math. Softw.", 22(4):469--483, 1996.

\bibitem{Berg2008} Berg M. de, Cheong O., Kreveld M. van, Overmars M.: \textit{"Computational Geometry: Algorithms and Applications"}. Springer (2008)

\bibitem{Branke2008} Branke J., Deb K., Miettinen K.: \textit{"Multiobjective Optimization: Interactive and Evolutionary Approaches"}. Springer (2008)

\bibitem{DLTZ} 
Deb K., Thiele L., Laumanns M., Zitzler E.: 
\textit{"Scalable multi-objective optimization test problems"}. 
" Evolutionary Computation", 1:825-830, 2002.

\bibitem{Ehrgott2005} Ehrgott M.: \textit{"Multicriteria Optimization"}. Springer (2005)

\bibitem{Ham} 
H\"am\"al\"ainen J.P., Miettinen K., Tarvainen P., Toivanen J.: 
\textit{"Interactive Solution Approach to a Multiobjective Optimization Problem in a Paper Machine Headbox Design"}. 
"Journal of Optimization Theory and Applications", 116:265, 2003.

\bibitem{Hapala2011} Hapala M. and Havran V.: \textit{"Review: Kd-tree Traversal Algorithms for Ray Tracing"}. "Computer Graphics Forum", 30(1):199-213, 2011.

\bibitem{Hartikainen2011} 
Hartikainen M., Miettinen K, Wiecek M.M.: 
\textit{"Constructing a Pareto front approximation for decision making"}. 
"Mathematical Methods of Operations Research", 73(2):209–234, 2011.

\bibitem{Hartikainen2012} 
Hartikainen M., Miettinen K, Wiecek M.M.: 
\textit{"PAINT: Pareto front interpolation for nonlinear multiobjective optimization"}. 
"Computational Optimization and Applications", 52(3):845-867, 2012.

\bibitem{Hartikainen2015} 
Hartikainen M. and Lovison A.: 
\textit{"PAINT-SiCon: constructing consistent parametric representations of Pareto sets in nonconvex multiobjective optimization"}. 
"Journal of Global Optimization", 62(2):243-261, 2015.

\bibitem{Harrington1965} 
Harrington, E.C.: 
\textit{"The Desirability Function"}. 
"Industrial Quality Control", 21:494-498, 1965.

\bibitem{Hwang1979} 
Hwang C.-L., Masud A.S.Md.: 
\textit{"Multiple Objective Decision Making — Methods and Applications: A State-of-the-Art Survey"}. 
"Mathematical Methods of Operations Research", Springer, Berlin, 1979.

\bibitem{Jar2004}
Jarre F. and Stoer J.: 
\textit{"Optimierung"}.
Springer, Berlin 2004.

\bibitem{NIMBUS} 
Miettinen K., M\"akel\"a M.M.: 
\textit{"Interactive bundle-based method for nondifferentiable multiobjeective optimization: NIMBUS"}. 
"Optimization", 34(3):231-146, 1995.

\bibitem{Miettinen1999} Miettinen K.: \textit{"Nonlinear Multiobjective Optimization"}. International series in operations researchand management science, Kluwer Academic Publishers (1999).

\bibitem{Monz2006} 
Monz M.: 
\textit{"Pareto Navigation
– interactive multiobjective
optimisation and its application
in radiotherapy planning"}. 
Ph.D. thesis, Technische Universität Kaiserslautern, 2006.

\bibitem{Monz2008} 
Monz M., K.H. K\"ufer, T.R. Bortfeld, C. Thieke: 
\textit{"Pareto navigation-algorithmic foundation of interactive multi-criteria IMRT planning"}. 
"Physics in Medicine and Biology", 53(4):985, 2008.

\bibitem{Pascoletti1984} 
Pascoletti A., Serafini P.: 
\textit{"Scalarizing vector optimization problems"}. 
"Journal of Optimization Theory and Applications", 42(4):499-524, 1984.

\bibitem{Ruzika2005} Ruzika S., Wiecek M. M.: \textit{"Approximation Methods in Multiobjective Programming"}. "Journal of Optimization Theory and Applications", 126(3):473-501, 2005.

\bibitem{Serna2008} 
Serna J. I.: 
\textit{"Approximating the Nondominated Set of R+ convex Bodies"}. 
Master thesis, Technische Universität Kaiserslautern, 2008.

\bibitem{Teichert2013} 
Teichert K.: 
\textit{"A hyperboxing Pareto approximation method applied to radiofrequency ablation treatment planning"}. 
Ph.D. thesis, Technische Universität Kaiserslautern, 2013.

\bibitem{VlENNET} 
R. VlENNET , C. FONTEIX, I. MARC, 
\textit{"Multicriteria optimization using a genetic algorithm for determining a Pareto set"}. 
"International Journal of Systems Science", 27(2):255-260, 1996.




\end{thebibliography}

\end{document}